\documentclass[11pt]{article}

\usepackage[reqno, namelimits, sumlimits]{amsmath}
\usepackage{amssymb, amsthm} 
\newtheorem{theorem}{Theorem} 
\newtheorem*{dlemma}{Doubling Lemma}
\newcommand{\R}{{\mathbb R}}
\newcommand{\eps}{\varepsilon}
\begin{document}

\title{Liouville theorems for scaling invariant superlinear
parabolic problems with gradient structure 
\thanks{Supported in part by
the Slovak Research and Development Agency under the contract No. APVV-0134-10 and by 
VEGA grant 1/0711/12.}
}
\author{Pavol Quittner \\ \\
 \small Department of Applied Mathematics and Statistics, Comenius University \\ 
 \small Mlynsk\'a dolina, 84248 Bratislava, Slovakia \\
 \tt quittner@fmph.uniba.sk  
}
\date{}

\maketitle

\begin{abstract} 
We provide a simple method for obtaining new Liouville theorems 
for scaling invariant superlinear parabolic problems with gradient structure.
To illustrate the method we prove Liouville theorems (guaranteeing
nonexistence of positive classical solutions) for the following model problems:
the scalar nonlinear heat equation
$$ u_t-\Delta u=u^p \qquad\hbox{in }\ \R^n\times\R, $$
its vector-valued generalization with a $p$-homogeneous nonlinearity  
and the linear heat equation in $\R^n_+\times\R$ complemented 
by nonlinear boundary conditions of the form $\partial u/\partial\nu=u^q$.
Here $\nu$ denotes the outer unit normal on the boundary of the halfspace $\R^n_+$
and the exponents $p,q>1$ satisfy 
$p<n/(n-2)$ and $q<(n-1)/(n-2)$ if $n>2$
(or $p<(n+2)/(n-2)$ and $q<n/(n-2)$ if $n>2$ and some symmetry of the solutions is assumed). 
As a typical application of our nonexistence results 
we provide optimal universal estimates for positive solutions
of related problems in bounded and unbounded domains.
\end{abstract}

\section{Introduction}
\label{intro}

In this paper we consider several  
model scaling invariant parabolic problems with gradient structure
and prove that these problems 
--- in a certain range of given parameters ---
do not possess positive entire solutions,
i.e. solutions defined for all times $t\in(-\infty,+\infty)$.
Such a result will be called (parabolic) Liouville theorem.

\paragraph{Nonlinear heat equation.} 
Let us first consider
the scalar nonlinear heat equation
\begin{equation} \label{eq-u}
 u_t-\Delta u=u^p, \qquad x\in\R^n,\ t\in\R, 
\end{equation}
where $p>1$, $n\geq1$ and $u=u(x,t)>0$.
Since problem (\ref{eq-u}) possesses positive stationary
solutions if $n>2$ and $p\geq(n+2)/(n-2)$,
the necessary condition for the Liouville theorem for (1) 
is $p<(n+2)/(n-2)_+$. This condition is also sufficient
if we restrict ourselves to radially symmetric solutions, see \cite{PQ,PQS}.
In the general non-radial case, the Liouville theorem for (\ref{eq-u}) was  
proved in \cite{BV}
only under the assumption $n=1$ or $n>1$ and $p<n(n+2)/(n-1)^2$.
In particular, if $n=2$ then one has to assume $p<8$.
Our main result for problem (\ref{eq-u})
guarantees that for $n=2$ this assumption on $p$ is superfluous.
More precisely, we prove the following Liouville theorem.
  
\begin{theorem} \label{thm1}
Let $p>1$, $(n-2)p<n$. Then the equation \hbox{\rm(\ref{eq-u})}
does not possess positive classical solutions.
\end{theorem}

If $n>2$ then $n/(n-2)<n(n+2)/(n-1)^2$ so that
the assertion in Theorem~\ref{thm1} follows
from \cite{BV} whenever $n\ne2$.
We formulate and prove our result for general $n$
since our method is very different from that in \cite{BV} and
it can also be used for problems where the arguments of \cite{BV}
cannot be used or have not been used so far.   
In particular, in this paper we also consider
a vector-valued generalization of (\ref{eq-u}) and
the linear heat equation complemented by nonlinear boundary conditions
and in these cases we obtain new results for all $n\geq1$.
It should be emphasized that we do not exploit the semilinear structure
of our problems: we consider these model problems just for simplicity. 
Our method is based on scaling and energy estimates for the rescaled problems.
This approach enables us to show that any positive bounded entire solution
of the parabolic problem has to be time-independent so that the nonexistence 
result for bounded solutions follows from the corresponding elliptic Liouville theorem
(and then the nonexistence of unbounded solutions is often an easy consequence
of doubling and scaling arguments). 
Let us note that if $n>2$ and $p>(n+2)/(n-2)$ then,
in addition to positive bounded stationary solutions, there also exist
positive bounded entire solutions of (\ref{eq-u}) which do depend on time;
in particular there exist homoclinic solutions, 
see \cite{FY}.

Liouville theorems  have important consequences concerning universal 
a priori estimates for positive solutions of related problems. 
To be more specific, let us formulate a typical result
of this type based on Theorem~\ref{thm1}. 
Since our result in Theorem~\ref{thm1} is new only if $n=2$, we restrict
ourselves to this case.
Consider nonnegative solutions of the equation
\begin{equation} \label{eq-f}
u_t-\Delta u=f(u) \qquad\hbox{in }\ \Omega\times(T_1,T_2),
\end{equation}
where $f:[0,\infty)\to\R$ is a continuous function satisfying
\begin{equation} \label{fp}
\lim_{u\to+\infty}u^{-p}f(u)=\ell\in(0,\infty),
\end{equation}
and
\begin{equation} \label{OmT}
\hbox{$\Omega$ is an arbitrary domain in $\R^2$, \ $-\infty\leq T_1<T_2\leq\infty$.}
\end{equation}
The following theorem is a direct consequence of Theorem~\ref{thm1} and
(the proof of) \cite[Theorems~3.1 and 4.1]{PQS}; cf.~also
\cite[Remark 3.4(e)]{PQS}.

\begin{theorem} \label{thm2}
Assume $p>1$, \hbox{\rm(\ref{fp}), (\ref{OmT})}
and let $u$ be a nonnegative classical solution of 
\hbox{\rm(\ref{eq-f})}. Then
\begin{equation} \label{ub}
 u(x,t)\leq C\bigl(C_1+(t-T_1)^{-\beta}+(T_2-t)^{-\beta}
  +C_2\hbox{\rm dist}^{-2\beta}(x,\partial\Omega)\bigr)
\quad\hbox{in }\ \Omega\times(T_1,T_2), 
\end{equation}
where $\beta:=1/(p-1)$, $C=C(f)>0$ is independent of $\Omega$, $T_1$, $T_2$ and $u$,
$C_1=0$ if $f(u)=u^p$, $C_1=1$ otherwise,
$C_2=1$, 
$(t-T_1)^{-\beta}:=0$ if $T_1=-\infty$,
$(T_2-t)^{-\beta}:=0$ if $T_2=\infty$ and
$\hbox{\rm dist}^{-2\beta}(x,\partial\Omega):=0$ if $\Omega=\R^2$.

If, in addition, $\Omega$ is (uniformly $C^2$) smooth
and $u$ satisfies the boundary condition 
\begin{equation} \label{Dbc}
u=0 \qquad{on }\ \partial\Omega\times(T_1,T_2)
\end{equation} 
then \hbox{\rm(\ref{ub})} is true with $C=C(f,\Omega)$, $C_1=1$ and $C_2=0$. 
\end{theorem}

In particular, if $\Omega\subset\R^2$ is smooth
and $u$ is any positive solution
of the problem (\ref{eq-f}),(\ref{Dbc}) which blows up at $t=T_2$
then Theorem~\ref{thm2} guarantees that the blow-up rate is of type I
and the corresponding estimate is universal 
(i.e. the constant $C$ in (\ref{ub}) does not depend on $u$).

Another application of Theorem~2 deals with so called ancient solutions.
Assume $T\in\R$, $1<p$ and $(n-2)p<n+2$.
Then \cite[Corollary 1.6]{MZ} gives a complete characterization of all
(positive classical) solutions of the problem
\begin{equation} \label{eq-uT}
 u_t-\Delta u=u^p, \qquad x\in\R^n,\ t\in(-\infty,T), 
\end{equation}
under the assumption
\begin{equation} \label{ass-MZ}
 u(x,t)\leq C(T-t)^{-\beta}.
\end{equation}   
Theorem~\ref{thm2} 
guarantees that (\ref{ass-MZ}) is always true if $n=2$.
In fact, the assertions in Theorem~\ref{thm2} 
(hence also (\ref{ass-MZ})) are true for any $n$ and $p>1$
such that \eqref{eq-u} does not possess positive classical solutions.

\paragraph{Vector valued case.}
Our next model problem is a vector-valued generalization of (\ref{eq-u}): 
we consider positive classical
solutions $U=(u_1,u_2,\dots,u_m)$ of the system
\begin{equation} \label{eq-U}
 U_t-\Delta U=F(U), \qquad x\in\R^n,\ t\in\R,
\end{equation}
where 
\begin{equation} \label{FG}
F=\nabla G,\ \hbox{ with }\ 
G\in C^{2+\alpha}_{loc}(\R^m,\R)\ \hbox{ for some }\ \alpha>0,
\end{equation} 
\begin{equation} \label{G}
G(0)=0,\qquad G(U)>0\quad\hbox{for }\ U\neq0,
\end{equation}
\begin{equation} \label{F1}
F(\lambda U)=\lambda^p F(U)\quad\hbox{for }\ U\in\R^m,\ \lambda>0,
\end{equation}
\begin{equation} \label{F2}
\hbox{there exists }\ 
\xi\in(0,\infty)^m \ \hbox{ such that }\ \ \xi\cdot F(U)>0\quad\hbox{for }\ U\ne0.
\end{equation}
Using the same arguments as in the proof of Theorem~1
we prove the following theorem.

\begin{theorem} \label{thmU}
Assume 
\hbox{\rm(\ref{FG})}, \hbox{\rm(\ref{G})}, 
\hbox{\rm(\ref{F1})}, \hbox{\rm(\ref{F2})} and
$p>1$, $(n-2)p<n$. Then the system \hbox{\rm(\ref{eq-U})}
does not possess nontrivial nonnegative classical solutions.
\end{theorem}

Notice that Theorem~\ref{thm1} is a special case of
Theorem~\ref{thmU}.
We will first prove Theorem~\ref{thm1}
(in order to explain the idea of our method by using the simplest possible
model problem);
the proof of Theorem~\ref{thmU} will then
follow the proof of Theorem~\ref{thm1}. 

Theorem~\ref{thmU} for $n=1$ and the approach in \cite[Proposition 2.4]{BPQ} 
(see also \cite{QS-D} and \cite{Phan})
enable us to prove also the following theorem.

\begin{theorem} \label{thmUa}
Assume 
\hbox{\rm(\ref{FG})}, \hbox{\rm(\ref{G})}, 
\hbox{\rm(\ref{F1})}, \hbox{\rm(\ref{F2})} and
$p>1$, $(n-2)p<n+2$. Then the system \hbox{\rm(\ref{eq-U})}    
does not possess nontrivial nonnegative classical radially symmetric solutions.
\end{theorem}

Theorem~\ref{thmUa} is a generalization of the scalar parabolic
Liouville theorem for radially symmetric solutions of (\ref{eq-u})
proved in \cite{PQ,PQS}
by completely different arguments.
Similarly as in the scalar case,
Theorems~\ref{thmU} and \ref{thmUa} can be used
in order to prove universal a priori estimates of positive
solutions of many related problems.

As far as we know, if $n,m>1$ then the only known nonexistence results 
for (\ref{eq-U}) in the non-radial case
are of Fujita-type and require the strong condition $p\leq (n+2)/n$.
If $n=1$, $m=2$ and
\begin{equation} \label{Fex} 
F(u_1,u_2)=(u_1^p-\lambda u_1^ru_2^{r+1},u_2^p-\lambda u_1^{r+1}u_2^r),
\qquad p=2r+1>1,
\end{equation}
then by using the approach in \cite{BV}, a Liouville theorem for
(\ref{eq-U}) has very recently been established in \cite{Phan}
under the assumption $\lambda<r/(3r+2)$.
Notice that in this particular case,
Theorem~\ref{thmU} guarantees the nonexistence
for any $\lambda<1$ and this condition on $\lambda$ is optimal.

In the radial setting, assuming $m=2$, (\ref{Fex})
and either $p=3\geq n$, $\lambda<1$ or $p(n-2)<n+2$, $\lambda<r/(3r+2)$, 
nonexistence results for (\ref{eq-U})
have also been obtained in \cite{QS-D} or \cite{Phan}, respectively.  

\paragraph{Nonlinear boundary conditions.}
Next consider positive classical solutions of the problem
\begin{equation} \label{eq-ubc}
\left.\begin{aligned}
 u_t-\Delta u &=0 &\qquad&\hbox{in }\R^n_+\times\R, \\
 u_\nu &= u^q      &\qquad&\hbox{on }\partial\R^n_+\times\R,
\end{aligned}\quad\right\}
\end{equation}
where  
$\R^n_+:=\{(x=(x_1,x_2,\dots,x_n)\in\R^n:x_1>0\}$, 
$\nu=(-1,0,0,\dots,0)$ is the outer unit normal on
the boundary $\partial\R^n_+=\{x\in\R^n:x_1=0\}$
and $q>1$.
In this case our method yields the following result.

\begin{theorem} \label{thm3}
Let $q>1$, $(n-2)q<n-1$. Then the problem \hbox{\rm(\ref{eq-ubc})}
does not possess positive classical bounded solutions.
\end{theorem}

The result in Theorem~\ref{thm3} is new for any $n\geq1$.
If $n=1$ then this nonexistence result was proved in \cite{QS-D} 
by other arguments, but only 
for solutions with bounded spatial derivatives.
For general $n\geq1$ the only known nonexistence results for (\ref{eq-ubc})
are of Fujita-type and require $q\leq(n+1)/n$,
see \cite{GL,DFL}.

Liouville theorem for stationary solutions
of (\ref{eq-ubc}) is true for $q<n/(n-2)_+$ (see \cite{Hu})
and this condition on $q$ is optimal: if $n>2$ and $q=n/(n-2)$
then there exists a stationary solution of (\ref{eq-ubc})
of the form $u(x)=c|x-x_0|^{2-n}$, where the first component
of $x_0$ is negative (see \cite{Har} and the references therein
for the analysis of stationary solutions for $q\geq n/(n-2)$).
Under the optimal assumption  $q<n/(n-2)_+$ we can also prove nonexistence
of solutions of (\ref{eq-ubc}) exhibiting the following axial symmetry:
\begin{equation} \label{sym}
u(x_1,\tilde x,t)=v(x_1,|\tilde x|,t), \quad
\hbox{where }\ \tilde x=(x_2,x_3,\dots,x_n).
\end{equation}

\begin{theorem} \label{thm3a}
Let $q>1$, $(n-2)q<n$. Then the problem \hbox{\rm(\ref{eq-ubc})}
does not possess positive classical bounded solutions
exhibiting the symmetry property
\hbox{\rm(\ref{sym})}.
\end{theorem}

Theorem~\ref{thm3a} is an analogue to Theorem~\ref{thmUa} 
and is proved by similar but technically more advanced arguments.

Let us also mention that the boundedness assumptions 
in Theorems~\ref{thm3} and \ref{thm3a}
still allow applications 
based on doubling and scaling arguments and yielding
a priori estimates for positive solutions of related problems.
In particular, 
Theorem~\ref{thm3} can be used to obtain blow-up rate
estimates for positive solutions of the problem 
\begin{equation} \label{eq-nbc}
\left.\begin{aligned}
 u_t-\Delta u &=0 &\qquad& x\in\Omega,\ t\in(0,T), \\
 u_\nu &= u^q      &\qquad& x\in\partial\Omega,\ t\in(0,T), 
\end{aligned}\quad\right\}
\end{equation}
where  $\nu$ denotes the outer unit normal on the boundary $\partial\Omega$.
More precisely, we will prove the following theorem.

\begin{theorem} \label{thm3b}
Assume that $\Omega\subset\R^n$ is bounded and smooth, $q>1$, 
$(n-2)q<n-1$. 
Assume also that $u$ is a positive classical solution of \hbox{\rm(\ref{eq-nbc})} 
which blows up at $t=T$. 
Then there exists $C=C(u)>0$ such that $u$
satisfies the blow-up rate estimate
\begin{equation} \label{est-rate}
u(x,t)(T-t)^{1/2(q-1)}+|\nabla u(x,t)|(T-t)^{q/2(q-1)}\leq C
\end{equation}
for all $x\in\overline\Omega$ and $t\in(T/2,T)$.
\end{theorem}

If $\Omega$ is bounded and $q>1$ then any positive solution of (\ref{eq-nbc})
blows up in finite time.
Theorem~\ref{thm3b} guarantees
that for $1<q<(n-1)/(n-2)_+$,
the blow-up of such solution 
is of type I, i.e.~$u$ satisfies the estimate 
\begin{equation} \label{rate-nbc}
\|u(\cdot,t)\|_\infty\leq C(T-t)^{-1/2(q-1)}\quad\hbox{for}\quad t\in(T/2,T).
\end{equation}
This result for bounded domains was known only under the stronger
assumption $1<q\leq1+1/n$ (see \cite{Hu96}).
On the other hand,
type I blow-up for both positive and sign-changing solutions of (\ref{eq-nbc}) has
been established in the full subcritical range $1<q<n/(n-2)_+$
if $\Omega$ is a half-space (see \cite{CF} and \cite{QS-P})
and it has also been proved for bounded domains and 
$1<q\leq n/(n-2)_+$
in the class of positive, time increasing solutions (see \cite{Hu96}).
Let us also mention that the blow-up rate estimate
(\ref{rate-nbc}) is optimal (see the lower estimates in \cite{HY,IS})
and that the blow-up need not be of type I 
for (some) supercritical $q$ (see \cite{Har-NA}).

\section{Proof of Theorems~\ref{thm1}, \ref{thmU} and \ref{thmUa}}  \label{sec-proof}

In the proofs we will often need the following lemma.

\begin{dlemma}
{\rm (see \cite[Lemma 5.1]{PQS}).}
Let $(X,d)$ be a complete metric space and
$\emptyset\ne D\subset X$.
Let $M:D\to(0,\infty)$ be bounded on compact subsets of $D$
and fix a real $k>0$. If  $y\in D$ is such that
\begin{equation} \label{M2k}
M(y)\,{\rm dist}(y,X\setminus D)>2k, 
\end{equation}
then there exists $x\in D$ such that
\begin{equation} \label{Mxy}
M(x)\,{\rm dist}(x,X\setminus D)>2k,\qquad M(x)\geq M(y), 
\end{equation}
and
\begin{equation} \label{M2M}
M(z) \leq 2M(x)\quad
\hbox{whenever }\ \hbox{\rm dist}(z,x)\leq \frac{k}{M(x)}.
\end{equation}
\end{dlemma}

Notice that the inequalities in (\ref{Mxy}) and (\ref{M2M})
guarantee
$$ \hbox{\rm dist}(z,x)\leq \frac{k}{M(x)}
 <\frac12\hbox{\rm dist}(x,X\setminus D),$$
so that $z\in D$ and the value $M(z)$ is well defined.
Notice also that
if $D=X$ then $\hbox{dist}(y,X\setminus D)=\hbox{dist}(y,\emptyset)=\infty$
so that the assumption (\ref{M2k}) is satisfied for any $y\in D$.
In most cases, we will use the Doubling Lemma 
with $X$ being a closed subset of $\R^n\times\R$
equipped with the parabolic distance
$$ \hbox{dist}_P((x,t),(\tilde x,\tilde t)):=
  |x-\tilde x|+\sqrt{|t-\tilde t|}.$$

\vskip5mm
{\it Proof of Theorem~\ref{thm1}.\/}
Assume on the contrary that there exists
a positive solution $u$ of (\ref{eq-u}).
Doubling and scaling arguments in \cite{PQS}
guarantee that we may assume that 
\begin{equation} \label{bound-u}
 u(x,t)\leq 1 \qquad\hbox{ for all }x\in\R^n,\ t\in\R.
\end{equation}
In fact, assume that
$u(x_0,t_0)>1$ for some $(x_0,t_0)\in\R^n\times\R$.
For any $k=1,2,\dots$, 
the Doubling Lemma 
(with $D=X=\R^n\times\R$,
$\hbox{dist}=\hbox{dist}_P$,
$M=u^{(p-1)/2}$
and $y=(x_0,t_0)$) guarantees the existence
of $(x_k,t_k)$ such that 
$$ \begin{aligned}
 M_k:=u^{(p-1)/2}(x_k,t_k) &\geq u^{(p-1)/2}(x_0,t_0)\quad\hbox{and}\\
  u^{(p-1)/2}(x,t) &\leq 2M_k \quad\hbox{whenever}\quad 
 |x-x_k|+\sqrt{|t-t_k|}\leq \frac{k}{M_k}.
\end{aligned}
$$
The rescaled functions
$$ v_k(y,s):=\lambda_k^{2/(p-1)} u(x_k+\lambda_k y,t_k+\lambda_k^2 s),
 \quad\hbox{where}\quad \lambda_k=\frac1{2M_k},$$
are positive solutions of (\ref{eq-u}) and satisfy $v_k(0,0)=2^{-2/(p-1)}$,
$v_k(y,s)\leq1$ for $|y|+\sqrt{|s|}\leq 2k$. 
The parabolic regularity guarantees that the sequence $\{v_k\}$
is relatively compact (in $C_{loc}$, for example), so that
a suitable subsequence of $\{v_k\}$
converges to a nonnegative solution $v$ of (\ref{eq-u}) satisfying $v\leq1$.
Since $v(0,0)>0$, we have $v>0$ by the maximum principle and uniqueness. 
Consequently, replacing $u$ by $v$ we may assume that (\ref{bound-u}) is true.

Denote $c_0:=u(0,0)$ and $\beta:=1/(p-1)$. For $y\in\R^n$, $s\in\R$ and
$k=1,2,\dots$ set 
$$w_k(y,s):=(k-t)^\beta u(y\sqrt{k-t},t),\qquad\hbox{where }\  s=-\log(k-t),\ \ t<k.$$
Set also $s_k:=-\log k$ and notice that
$w=w_k$ solve the problem
\begin{equation} \label{eq-w}
\left.\begin{aligned}
 w_s &=\Delta w-\frac12 y\cdot\nabla w-\beta w+w^p \\
     &=\frac1\rho\nabla\cdot(\rho\nabla w)-\beta w+w^p\qquad \hbox{in }
\R^n\times\R,
\end{aligned}\quad\right\}
\end{equation}
where $\rho(y):=e^{-|y|^2/4}$. In addition,
$$  w_k(0,s_k)=k^\beta c_0 $$  
and
\begin{equation} \label{bound-w2}
 \|w_k(\cdot,s)\|_\infty\leq e^{2\beta}k^\beta \qquad \hbox{for }\
s\in[s_k-2,\infty).
\end{equation}
Set
$$ E_k(s):=\frac12\int_{\R^n}(|\nabla w_k(y,s)|^2+\beta w_k^2(y,s))\rho(y)\,dy
 -\frac1{p+1}\int_{\R^n}w_k^{p+1}(y,s)\rho(y)\,dy.$$
Then in the same way as in \cite[(2.25) and Proposition~2.1]{GK}
we obtain $E_k(s)\geq0$ and, given $\sigma<s_k$,
\begin{equation} \label{GK1}
\left.\begin{aligned}
\frac12&\Bigl(\int_{\R^n}w_k^2(y,s_k)\rho(y)\,dy
             -\int_{\R^n}w_k^2(y,\sigma)\rho(y)\,dy\Bigr) \\
 &= -2\int_{\sigma}^{s_k}E_k(s)\,ds
  +\frac{p-1}{p+1}\int_{\sigma}^{s_k}\int_{\R^n}
    w_k^{p+1}(y,s)\rho(y)\,dy\,ds,
\end{aligned} \quad\right\}
\end{equation}
\begin{equation} \label{GK2}
\int_{\sigma}^{s_k}\int_{\R^n} 
\Big|\frac{\partial w_k}{\partial s}(y,s)\Big|^2\rho(y)\,dy\,ds
= E_k(\sigma)-E_k(s_k)\leq E_k(\sigma).
\end{equation}
Multiplying equation (\ref{eq-w}) by $\rho$, integrating over
$y\in\R^n$ and using Jensen's inequality yields
$$ \begin{aligned}
\frac{d}{ds} &\int_{\R^n}w_k(y,s)\rho(y)\,dy+ \beta \int_{\R^n}w_k(y,s)\rho(y)\,dy \\
 &=\int_{\R^n}w_k^p(y,s)\rho(y)\,dy
 \geq C_{n,p}\Bigl(\int_{\R^n}w_k(y,s)\rho(y)\,dy\Bigr)^p,
\end{aligned}
 $$
where $C_{n,p}:=(4\pi)^{-n(p-1)/2}$,
which (as in the proof of \cite[Theorem 1]{FSW}, for example)
implies the estimates
\begin{equation} \label{GK4}
\int_{\R^n}w_k(y,s)\rho(y)\,dy\leq \tilde C_{n,p}
\end{equation}
and
\begin{equation} \label{GK5}
\int_{\sigma}^{s_k}\int_{\R^n}w_k^p(y,s)\rho(y)\,dy\,ds\leq
\tilde C_{n,p}(1+\beta(s_k-\sigma)),
\end{equation}
where $\tilde C_{n,p}=(\beta/C_{n,p})^\beta$.
The monotonicity  
of $E_k$,
(\ref{GK1}), (\ref{bound-w2}), (\ref{GK4}) and (\ref{GK5})
guarantee
$$ \begin{aligned}
2&E_k(s_k-1) 
\leq 2\int_{s_k-2}^{s_k-1} E_k(s)\,ds
\leq 2\int_{s_k-2}^{s_k} E_k(s)\,ds \\
 &\leq \frac12\int_{\R^n}w_k^2(y,s_k-2)\rho(y)\,dy
 +\frac{p-1}{p+1}\int_{s_k-2}^{s_k}\int_{\R^n}
   w_k^{p+1}(y,s)\rho(y)\,dy\,ds  \\
 &\leq e^{2\beta}k^\beta\Bigl( \int_{\R^n}w_k(y,s_k-2)\rho(y)\,dy
 +\int_{s_k-2}^{s_k}\int_{\R^n}
   w_k^{p}(y,s)\rho(y)\,dy\,ds\Bigr)  \\
&\leq 2C(n,p)k^\beta,
\end{aligned} $$ 
where $C(n,p):=e^{2\beta}\tilde C_{n,p}(1+\beta)$,
hence 
$E_k(s_k-1)\leq C(n,p)k^\beta$.
This estimate and (\ref{GK2}) guarantee
\begin{equation} \label{GK6}
\int_{s_k-1}^{s_k}\int_{\R^n} 
\Big|\frac{\partial w_k}{\partial s}(y,s)\Big|^2\rho(y)\,dy\,ds
\leq C(n,p)k^\beta.
\end{equation}
Denote $\lambda_k:=k^{-1/2}$ and set
$$ v_k(z,\tau):=\lambda_k^{2/(p-1)}w_k(\lambda_k z,\lambda_k^2\tau+s_k),
 \qquad z\in\R^n,\ -k\leq\tau\leq0.  $$ 
Then
$0<v_k\leq e^{2\beta}$, $v_k(0,0)=c_0$,
$$ \frac{\partial v_k}{\partial\tau}-\Delta v_k-v_k^p
  =-\lambda_k^2\Bigl(\frac12 z\cdot\nabla v_k+\beta v_k\Bigr) $$
and, denoting $\alpha:=-n+2+4/(p-1)$ and using (\ref{GK6}) we also have
$$
\begin{aligned}
\int_{-k}^0\int_{|z|<\sqrt{k}}
 \Big|\frac{\partial v_k}{\partial\tau}(z,\tau)\Big|^2\,dz\,d\tau
 &=\lambda_k^\alpha
 \int_{s_k-1}^{s_k}\int_{|y|<1} 
\Big|\frac{\partial w_k}{\partial s}(y,s)\Big|^2\,dy\,ds \\
&\leq C(n,p)e^{1/4}k^{-\alpha/2+\beta}\to 0 \quad\hbox{as }\
k\to\infty.
\end{aligned}$$
Now the same arguments as in \cite{GK} show that
(up to a subsequence) the sequence $\{v_k\}$
converges to a positive solution $v=v(z)$
of the problem $\Delta v+v^p=0$ in $\R^n$,
which contradicts the elliptic Liouville theorem in \cite{GS}. 

Notice that the explicit formula
$$
v_k(z,\tau)=e^{-\beta\tau/k}u(e^{-\tau/2k}z,k(1-e^{-\tau/k}))
$$
guarantees $v_k\to u$.
Notice also that if $p=n/(n-2)$ and if
we rescaled the functions $w_k$ on the intervals
$[s_k-1,s_k+1]$ instead of $[s_k-1,s_k]$ then the above arguments would
guarantee $\int_\infty^\infty\int_{\R^n}u_t^2\,dx\,dt\leq C(n,p)e^{1/4}$.
\qed

\vskip5mm
{\it Proof of Theorem~\ref{thmU}.\/}
Assume on the contrary that there exists
a nontrivial nonnegative solution $U$ of (\ref{eq-U}).
As in the proof of Theorem~\ref{thm1},
doubling and scaling arguments in \cite{PQS}
guarantee that we may assume 
$$ 
 |U(x,t)|\leq 1 \qquad\hbox{ for all }x\in\R^n,\ t\in\R.
$$ 
Denote $C_0:=U(0,0)$ and $\beta:=1/(p-1)$. For $y\in\R^n$, $s\in\R$ and
$k=1,2,\dots$ set 
$$W_k(y,s):=(k-t)^\beta U(y\sqrt{k-t},t),\qquad\hbox{where }\  s=-\log(k-t),\ \ t<k.$$
Set also $s_k:=-\log k$ and notice that 
$W=W_k$ solve the problem
\begin{equation} \label{eq-W}
\left.\begin{aligned}
 W_s &=\Delta W-\frac12 y\cdot\nabla W-\beta W+F(W) \\
     &=\frac1\rho\nabla\cdot(\rho\nabla W)-\beta W+F(W)\qquad \hbox{in }
\R^n\times\R,
\end{aligned}\quad\right\}
\end{equation}
where $\rho(y)=e^{-|y|^2/4}$. In addition,
$$ 
 W_k(0,s_k)=k^\beta C_0
\quad\hbox{and}\quad
 \|W_k(\cdot,s)\|_\infty\leq e^{2\beta}k^\beta \qquad \hbox{for }\
s\in[s_k-2,\infty).
$$
Set
$$ E_k(s):=\frac12\int_{\R^n}
 (|\nabla W(y,s)|^2+\beta W^2(y,s))\rho(y)\,dy
 -\int_{\R^n}G(W(y,s))\rho(y)\,dy.$$
Since assumptions 
\hbox{\rm(\ref{FG})},
(\ref{G}), (\ref{F1}) and (\ref{F2}) guarantee
$$C_G|W|^{p+1}\geq G(W)=\frac1{p+1}F(W)W\geq c_G|W|^{p+1},\qquad
\xi\cdot F(W)\geq c_F|W|^p,$$ 
one can use 
the same arguments as in the proof of Theorem~\ref{thm1}
to show that $E_k(s_k-1)\leq Ck^\beta$
for some $C$ depending only on $n,p,C_G,c_F$ and $\xi$. 
In fact, to prove the analogoues of (\ref{GK4}) and (\ref{GK5}), for example,
it is sufficient to multiply the $i$-th component in (\ref{eq-W})
by $\xi_i\rho$, integrate and sum over $i$.
Consequently, as in the proof of Theorem~\ref{thm1} the functions 
$$ V_k(z,\tau):=\lambda_k^{2/(p-1)}W_k(\lambda_k z,\lambda_k^2\tau+s_k),
 \qquad z\in\R^n,\ -k\leq\tau\leq0  $$
converge (up to a subsequence) to a positive solution $V=V(z)$
of the problem $\Delta V+F(W)=0$ in $\R^n$,
which contradicts the elliptic Liouville theorem  \cite[Theorem 6(i)]{QS9}. 
\qed

\vskip5mm
{\it Proof of Theorem~\ref{thmUa}.\/}
The proof is based on the same arguments as the proof of
\cite[Proposition 2.4]{BPQ} 
(cf.~also \cite[Theorem 4.1]{QS-D}).
For the reader's convenience (and since we will also need
a nontrivial modification of these arguments in the
proof of Theorem~\ref{thm3a}) we provide a detailed proof.

Let $U$ be a nontrivial nonnegative radial solution of (\ref{eq-U}).
Since $U$  is radial, there exists 
$\tilde U:[0,\infty)\times\R\to\R^m:(r,t)\mapsto\tilde U(r,t)$
such that $U(x,t)=\tilde U(|x|,t)$.

First we show that we can assume that $U$ is bounded.
In fact, assume that there exist $r_k\in[0,\infty)$ and $t_k\in\R$
such that $|\tilde U(r_k,t_k)|\to\infty$. 
The Doubling Lemma 
(with $D=X=[0,\infty)\times\R$, 
$\hbox{dist}=\hbox{dist}_P$
and $M=|\tilde U|^{(p-1)/2}$)
guarantees
that we may assume
$$ M(r,t)\leq 2M_k
\quad\hbox{whenever}\quad 
|r-r_k|+\sqrt{|t-t_k|}\leq \frac{k}{M_k},$$
where $M_k:=|\tilde U(r_k,t_k)|^{(p-1)/2}$.
Set $\rho_k:=r_kM_k$ and $\lambda_k:=1/M_k$.
Passing to a subsequence we may assume
$\rho_k\to\rho_\infty\in[0,\infty]$.
If $\rho_\infty=\infty$ then
the functions
$$V_k(\rho,s):=
\lambda_k^{2/(p-1)}
 \tilde U(r_k+\lambda_k\rho,t_k+\lambda_k^2s), \qquad \rho\geq-\rho_k,\ s\in\R, 
$$
solve the equations
$$ \partial_t V_k-\partial_{\rho\rho}V_k=
 \frac{n-1}{\rho_k+\rho}\partial_\rho V_k +F(V_k)$$ 
and a subsequence of $\{V_k\}$ converges to nontrivial nonnegative
solution of (\ref{eq-U}) with $n=1$, which contradicts Theorem~\ref{thmU}.
Hence $\rho_\infty<\infty$.
The functions
$$V_k(\rho,s):=
\lambda_k^{2/(p-1)}
 \tilde U(\lambda_k\rho,t_k+\lambda_k^2s), 
\qquad
\rho\geq0,\ s\in\R,
$$
solve the equations
$$ \partial_t V_k-\partial_{\rho\rho}V_k=
 \frac{n-1}{\rho}\partial_\rho V_k +F(V_k)
$$
and satisfy $|V_k(\rho_k,0)|=1$,
$|V_k(\rho,s)|\leq 2^{2/(p-1)}$ for $|\rho-\rho_k|+\sqrt{|s|}\leq k$.
Passing to a subsequence
we may assume $V_k\to V$, where $V$ is a nontrivial nonnegative bounded radial 
solution of of (\ref{eq-U}).
Replacing $U$ by $V$ we may assume that $U$ is bounded.

Since $U$ is bounded, 
the parabolic regularity implies that $\nabla U$ is bounded as well, hence
\begin{equation} \label{bound-U}
 |U|+|\nabla U|\leq C.
\end{equation}
Now we use similar doubling and scaling arguments as above
to show the uniform decay estimate
\begin{equation} \label{decayU}  
 |\tilde U(r,t)|r^{2/(p-1)}+|\nabla\tilde U(r,t)|r^{(p+1)/(p-1)}\leq C
\end{equation}
(where the constant $C$ is different from that in (\ref{bound-U})).
Assume on the contrary that there exist $r_k>0$ and $t_k\in\R$
such that
$$ |\tilde U(r_k,t_k)|r_k^{2/(p-1)}+|\nabla\tilde U_r(r_k,t_k)|r_k^{(p+1)/(p-1)}\to\infty.$$
Set
$$ M(r,t):=|\tilde U(r,t)|^{(p-1)/2}+|\nabla\tilde U_r(r,t)|^{(p-1)/(p+1)},
 \qquad r>0,\ t\in\R$$
and $M_k:=M(r_k,t_k)$.
Then $M_kr_k\to\infty$ and
passing to a subsequence we may assume
$M_k> 2k/r_k$.              
The Doubling Lemma 
(with $X=[0,\infty)\times\R$, $D=(0,\infty)\times\R$ and
$\hbox{dist}=\hbox{dist}_P$)
guarantees that we may assume
$$ M(r,t)\leq 2M_k \qquad\hbox{whenever}\quad
|r-r_k|+\sqrt{|t-t_k|}\leq\frac k{M_k}. $$
Set $\lambda_k:=1/M_k$ and
$$V_k(\rho,s):=\lambda_k^{2/(p-1)}
 \tilde U(r_k+\lambda_k\rho,t_k+\lambda_k^2s).$$
Then 
$$\begin{aligned}
|V_k(0,0)|^{(p-1)/2}+|\partial_\rho V_k(0,0)|^{(p-1)/(p+1)}&=1, \\
|V_k(\rho,s)|^{(p-1)/2}+|\partial_\rho V_k(\rho,s)|^{(p-1)/(p+1)}&\leq 2
\quad
\hbox{whenever }\ |\rho|+\sqrt{|s|}\leq k,
\end{aligned}$$
and $V_k$ solves the equation
$$ \partial_t V_k-\partial_{\rho\rho}V_k=
 \frac{n-1}{r_k/\lambda_k+\rho}\partial_\rho V_k +F(V_k).$$ 
Since $r_k/\lambda_k=r_kM_k\to\infty$, 
it is easy to pass to the limit to get a nontrivial nonnegative bounded solution $V$
of \eqref{eq-U} with $n=1$.
However, this contradicts Theorem~\ref{thmU}.
Consequently, (\ref{decayU}) is true.

Next we use the energy functional
$$ E(U(\cdot,t)):=\int_{\R^n}\Bigl(\frac12|\nabla U(x,t)|^2-G(U(x,t))\Bigr)\,dx. $$
The arguments in \cite[Example 51.28, the case $\lambda=0$]{QS})
guarantee that the system (\ref{eq-U}) is well posed in
the space
\begin{equation}
  \label{e-space}
\left.\begin{aligned}
  {\cal E} &:=\{W\in  L^{p+1}(\R^n,\R^m): \nabla W\in L^2(\R^n,\R^{mn})\},\\
 \|W\|_{\cal E} &:=\|W\|_{L^{p+1}}+\|\nabla W\|_{L^2}
\end{aligned}\qquad\right\}
\end{equation}
and the corresponding solution satisfies the energy identity  
\begin{equation}
  \label{e-identity}
  E(U(\cdot,t_2))-E(U(\cdot,t_1))
        = -\int_{t_1}^{t_2}\int_{\R^n}|U_t|^2(x,t)\,dx\,dt.
\end{equation}
Estimates (\ref{decayU}) and (\ref{bound-U}) guarantee $\|U(\cdot,t)\|_{\cal E}\leq C$
and $|E(u(\cdot,t))|\leq C$ with $C$ independent of $t$, hence
$$\int_{\R}\int_{\R^n}|U_t|^2\,dx\,dt<\infty$$
and
\begin{equation} \label{utk}
\int_{|t|>k}\int_{\R^n}|U_t|^2\,dx\,dt\to0
\quad\hbox{ as }\ k\to\infty.
\end{equation}
Next we claim
\begin{equation} \label{Usup}
 \sup_{x\in\R^n,\ |t|>2k}(|U(x,t)|+|\nabla U(x,t)|)\to0
\quad\hbox{ as }\ k\to\infty.
\end{equation}
Assume on the contrary that there exist $x_k\in\R^n$
and $t_k\in\R$, $|t_k|>2k$, such that
$$ |U(x_k,t_k)|+|\nabla U(x_k,t_k)|\geq c_0>0.$$
Estimate (\ref{decayU}) shows that the sequence $\{x_k\}$ is bounded
so that we may assume $x_k\to x_\infty$.
Set $V_k(x,t):=U(x,t-t_k)$. Then a subsequence of $\{V_k\}$ converges
(locally uniformly in $C^1$) 
to a nonnegative radial solution $V$ of \eqref{eq-U}.
Estimate $|V(x_\infty,0)|+|\nabla V(x_\infty,0)|\geq c_0$
shows that
$V$ is nontrivial
and estimate (\ref{utk}) guarantees that $V$ does not depend on $t$.
However, this contradicts the elliptic Liouville theorem
\cite[Proposition~5(i)]{QS9}.

Estimates (\ref{decayU}) and (\ref{Usup})
guarantee 
$E(U(\cdot,t))\to0$ as $|t|\to\infty$,
so that  $E(U(\cdot,t))\equiv0$ 
by the monotonicity of $t\mapsto E(U(\cdot,t))$.
Consequently, $U_t\equiv0$ which contradicts 
\cite[Proposition~5(i)]{QS9}.
\qed

\section{Proofs of Theorems~\ref{thm3}, \ref{thm3a} and \ref{thm3b}}  \label{sec-proof2}

{\it Proof of Theorem~\ref{thm3}.\/}
The proof will follow that of Theorem~\ref{thm1}
but we will also need some additional arguments.

Assume on the contrary that there exists
a positive bounded solution $u$ of (\ref{eq-ubc}).
By using doubling and scaling arguments we first show that 
we may assume 
\begin{equation} \label{bound-ubc2} 
 u(x,t)+|\nabla u(x,t)|\leq C \qquad\hbox{ for all }x\in\overline{\R^n_+},\ t\in\R.
\end{equation}
Assume that (\ref{bound-ubc2}) fails.
Since $u\leq C_u$ for some $C_u>0$,
we can find $(x_k,t_k)$ such that
$|\nabla u(x_k,t_k)|\to\infty$.
Set 
$$ M(x,t):=u^{q-1}(x,t)+|\nabla u(x,t)|^{(q-1)/q}, $$
$M_k:=M(x_k,t_k)$ and $\lambda_k:=1/M_k$.
The Doubling Lemma (with $X=D=\overline{\R^n_+}\times\R$
and $\hbox{dist}=\hbox{dist}_P$)
guarantees that we may assume
$$ M(x,t)\leq 2M_k\quad\hbox{whenever}\quad
   |x-x_k|+\sqrt{|t-t_k|}\leq \frac{k}{M_k}.$$
Passing to a subsequence we may assume $c_k:=x_{k,1}M_k\to c_\infty\in[0,\infty]$,
where $x_{k,1}$ denotes the first component of $x_k$.
If $c_\infty=\infty$ then setting
$$ v_k(y,s):=\lambda_k^{1/(q-1)}u(x_k+\lambda_k y,
         t_k+\lambda_k^2 s),
   \quad y\in\R^n,\ y_1\geq-c_k,\ s\in\R,$$
a suitable subsequence of $\{v_k\}$ converges
to a nonnegative bounded solution $v$ 
of the linear heat equation in $\R^n\times\R$ 
satisfying $|\nabla v(0,0)|=1$, which contradicts
the Liouville theorem for the linear heat equation
(see \cite[Theorem 1]{Eid} or \cite[Theorem 4]{Eid1}
and cf.~also \cite{Nic}).
Therefore we have $c_\infty<\infty$.
Set $x_k^0:=(0,x_{k,2},\dots,x_{k,n})$,
$y_k:=(c_k,0,0,\dots,0)$
and 
$$ v_k(y,s):=\lambda_k^{1/(q-1)}u(x_k^0+\lambda_k y,
         t_k+\lambda_k^2 s),
   \quad y\in\R^n_+,\ \ s\in\R.$$
Then $v_k^{q-1}(y_k,0)+|\nabla v_k(y_k,0)|^{(q-1)/q}=1$
and a suitable subsequence of $\{v_k\}$ converges
to a nonnegative nontrivial (hence positive) bounded solution $v$
of (\ref{eq-ubc})
with bounded spatial derivatives. Replacing $u$ by $v$
we obtain (\ref{bound-ubc2}). 

If (\ref{bound-ubc2}) is true then the function
$$ v(y,s)=\lambda^{1/(q-1)}u(\lambda y,\lambda^2 s)
  \quad\hbox{where}\quad \lambda^{1/(q-1)}=1/C $$
is a positive solution of (\ref{eq-ubc}) satisfying (\ref{bound-ubc2}) with
$C=1$. Hence, replacing $u$ by $v$ 
we may assume
\begin{equation} \label{bound-ubc1}
 u(x,t)+|\nabla u(x,t)|\leq 1 \qquad\hbox{ for all }x\in\overline{\R^n_+},\ t\in\R.
\end{equation}

Next we prove that
\begin{equation} \label{monotone}
 u_{x_1}(x,t)\leq0 \qquad\hbox{ for all }x\in\R^n_+,\ t\in\R.
\end{equation}
The function $z:=u_{x_1}$ is bounded and satisfies
$$ z_t-\Delta z=0 \quad\hbox{in }\R^n_+\times\R, \qquad
  z<0 \quad\hbox{on }\partial\R^n_+\times\R.$$
In order to prove (\ref{monotone}) it is sufficient to show
$z(x,t)\leq\eps x_1$ for any $\eps>0$. 
Fix $\eps>0$ and set $v(x,t):=z(x,t)-\eps x_1$. 
Since $z$ is bounded, there exists $\lambda=\lambda(\eps)>0$ such that
$v(x,t)<0$ for $x_1\geq\lambda$. 
To show $v(x,t)\leq0$ for $x_1<\lambda$ we will proceed similarly
as in the proof of \cite[Theorem 2.4]{PQS}.

Denoting
$T_\lambda:=\{x\in\R^n:0<x_1<\lambda\}$ the function $v$ satisfies
$$ v_t-\Delta v=0 \quad\hbox{in }T_\lambda\times\R, \qquad
  v(x,t)<0 \quad\hbox{on }\partial T_\lambda\times\R.$$
Choosing $q\in(0,\pi^2/\lambda^2)$, \cite{Dan} guarantees the existence
of a smooth positive function $h$ on $\overline{T_\lambda}$ such that
$$ \Delta h+qh=0 \quad\hbox{in }T_\lambda \qquad\hbox{and}\qquad
 h(x)\to+\infty \quad\hbox{as }|x|\to\infty,\ x\in\overline{T_\lambda}.$$
In particular $h(x)\geq h_0>0$.
Set $w:=e^{qt}v/h$. Then $w$ satisfies 
$$ w_t-\Delta w-\frac{2\nabla h}{h}\cdot\nabla w=0 \quad\hbox{in }T_\lambda\times\R,
\qquad w\leq 0 \quad\hbox{on }\partial T_\lambda\times\R. $$
Fix $t_0<t_1$ and consider $(x,t)\in \overline{T_\lambda}\times[t_0,t_1]$.
Then $w(x,t)\to0$ as $|x|\to\infty$ 
and the maximum principle guarantees 
$$  \sup_{x\in T_\lambda} w^-(x,t_1)\le  \sup_{x\in T_\lambda} w^-(x,t_0), $$
where $w^-(x,t):=-\min(w(x,t),0)$. For
$v$ the above inequality  means
$$   \sup_{x\in T_\lambda} \frac{v^-(x,t_1)}{h(x)}
  \le e^{-q(t_1-t_0)}  \sup_{x\in T_\lambda} \frac{v^-(x,t_0)}{h(x)}. $$
In view of boundedness of $v$ on $T_\lambda\times\R$, letting $t_0\to -\infty$ we
obtain that $v(x,t_1)\le 0$.
This concludes the proof of (\ref{monotone}).

Denote $c_0:=u(0,0)$ and $\tilde\beta:=1/2(q-1)$. For $y\in\R^n_+$, $s\in\R$ and
$k=1,2,\dots$ set 
$$w_k(y,s):=(k-t)^{\tilde\beta} u(y\sqrt{k-t},t),\qquad\hbox{where }\  s=-\log(k-t),\ \ t<k.$$
Set also $s_k:=-\log k$ and notice that
$w=w_k$ solve the problem
\begin{equation} \label{eq-wbc}
\left.\begin{aligned}
 w_s &=\Delta w-\frac12 y\cdot\nabla w-\tilde\beta w 
   =\frac1\rho\nabla\cdot(\rho\nabla w)-\tilde\beta w &\quad& \hbox{in }
\R^n_+\times\R, \\
 w_\nu &= w^q &\quad& \hbox{on }\partial\R^n_+\times\R, 
\end{aligned}\ \right\}
\end{equation}
where $\rho(y)=e^{-|y|^2/4}$. 
In addition,
$$ w_k(0,s_k)=k^{\tilde\beta} c_0, \qquad
 \|w_k(\cdot,s)\|_\infty\leq e^{2\tilde\beta}k^{\tilde\beta} \quad \hbox{for }\
s\in[s_k-2,\infty).
$$
Set
$$ E_k(s):=\frac12\int_{\R^n_+}(|\nabla w_k(y,s)|^2+\tilde\beta w_k^2(y,s))\rho(y)\,dy
 -\frac1{q+1}\int_{\partial\R^n_+}w_k^{q+1}(\xi,s)\rho(\xi)\,dS_\xi.$$
Then $E_k(s)\geq0$ (see \cite{CF})
and, given $\sigma<s_k$, we also have 
$$ 
\begin{aligned}
\frac12&\Bigl(\int_{\R^n_+}w_k^2(y,s_k)\rho(y)\,dy
             -\int_{\R^n_+}w_k^2(y,\sigma)\rho(y)\,dy\Bigr) \\
 &= -2\int_{\sigma}^{s_k}E_k(s)\,ds
  +\frac{q-1}{q+1}\int_{\sigma}^{s_k}\int_{\partial\R^n_+}
    w_k^{q+1}(\xi,s)\rho(\xi)\,dS_\xi\,ds,
\end{aligned}
$$ 
$$
\int_{\sigma}^{s_k}\int_{\R^n_+} 
\Big|\frac{\partial w_k}{\partial s}(y,s)\Big|^2\rho(y)\,dy\,ds
= E_k(\sigma)-E_k(s_k)\leq E_k(\sigma).
$$ 
Since (\ref{monotone}) guarantees $\partial w_k/\partial{y_1}\leq0$,
we have
$$\begin{aligned}
\sqrt{\pi/2}\int_{\partial\R^n_+}w_k(\xi,s)\rho(\xi)\,dS_\xi
 &= \int_{\R^n_+}w_k((0,y_2,y_3,\dots,y_n),s)\rho(y)\,dy \\
 &\geq \int_{\R^n_+}w_k(y,s)\rho(y)\,dy.
\end{aligned}$$
Consequently,
multiplying the equation in (\ref{eq-wbc}) by $\rho$ and integrating over
$y\in\R^n_+$ yields
$$ \begin{aligned}
\frac{d}{ds}\int_{\R^n_+} &w_k(y,s)\rho(y)\,dy
 +\tilde\beta \int_{\R^n_+}w_k(y,s)\rho(y)\,dy
 =\int_{\partial\R^n_+}w_k^q(\xi,s)\rho(\xi)\,dS_\xi\\
&\geq C_{n-1,q}\Bigl(\int_{\partial\R^n_+}w_k(\xi,s)\rho(\xi)\,dS_\xi\Bigr)^q 
\geq \hat C_{n-1,q}\Bigl(\int_{\R^n_+}w_k(y,s)\rho(y)\,dy\Bigr)^q,
\end{aligned} $$
which again implies the estimates of the type
$$ 
\begin{aligned}
\int_{\R^n_+}w_k(y,s)\rho(y)\,dy &\leq\tilde C_{n-1,q}, \\
\int_{\sigma}^{s_k}\int_{\partial\R^n_+}w_k^q(\xi,s)\rho(\xi)\,dS_\xi\,ds
&\leq \tilde C_{n-1,q}(1+\tilde\beta(s_k-\sigma)).
\end{aligned} 
$$ 
In the same way as in the proof of Theorem~\ref{thm1},
the estimates above guarantee
$E_k(s_k-1)\leq \tilde C(n-1,q)k^{\tilde\beta}$ for suitable $\tilde C(n-1,q)$
and, consequently,
\begin{equation} \label{GK6bc}
\int_{s_k-1}^{s_k}\int_{\R^n_+} 
\Big|\frac{\partial w_k}{\partial s}(y,s)\Big|^2\rho(y)\,dy\,ds
\leq \tilde C(n-1,q)k^{\tilde\beta}.
\end{equation}
Denote $\lambda_k:=k^{-1/2}$ and set
$$ v_k(z,\tau):=\lambda_k^{1/(q-1)}w_k(\lambda_k z,\lambda_k^2\tau+s_k),
 \qquad z\in\R^n_+,\ -k\leq\tau\leq0.  $$ 
Then
$0<v_k\leq e^{2\tilde\beta}$, $v_k(0,0)=c_0$,
$$ 
\begin{aligned}
\frac{\partial v_k}{\partial\tau}-\Delta v_k 
  &=-\lambda_k^2\Bigl(\frac12 z\cdot\nabla v_k+\tilde\beta v_k\Bigr)
&\quad&\hbox{in }\R^n_+\times(-k,0), \\
v_\nu &=v^q &\quad&\hbox{on }\partial\R^n_+\times(-k,0),
\end{aligned}
 $$
and, denoting $\tilde\alpha:=-n+2+2/(q-1)$ and using (\ref{GK6bc}) we also have
$$
\begin{aligned}
\int_{-k}^0\int_{|z|<\sqrt{k}}
 \Big|\frac{\partial v_k}{\partial\tau}(z,\tau)\Big|^2\,dz\,d\tau
 &=\lambda_k^{\tilde\alpha}
 \int_{s_k-1}^{s_k}\int_{|y|<1} 
\Big|\frac{\partial w_k}{\partial s}(y,s)\Big|^2\,dy\,ds \\
&\leq C(q)e^{1/4}k^{-\tilde\alpha/2+\tilde\beta}\to 0 \quad\hbox{as }\
k\to\infty.
\end{aligned}$$
As in the proof of Theorem~\ref{thm1} (cf.~also \cite{CF}),
a subsequence of $\{v_k\}$
converges to a positive solution $v=v(z)$
of the problem $\Delta v=0$ in $\R^n_+$, $v_\nu=v^q$ on $\partial\R^n_+$,
which contradicts the elliptic Liouville theorem in \cite{Hu}. 
\qed

\vskip5mm
{\it Proof of Theorem~\ref{thm3a}.\/}
Due to Theorem~\ref{thm3} we may assume $n>2$ and $n-1\leq q(n-2)<n$.
Assume that $u$ is a positive classical bounded solution of (\ref{eq-ubc})
satisfying (\ref{sym}).
Similarly as in the proof of Theorem~\ref{thm3}
we will first show that we may assume 
\begin{equation}  \label{bound-unu}
 u+|\nabla u|\leq C
\end{equation}
and then (similarly as in the proof of Theorem~\ref{thmUa})
we will prove that $u$ satisfies
suitable decay estimates which allow us to use the energy functional
\begin{equation} \label{E-u}
E(\varphi):=\frac12\int_{\R^n_+}|\nabla\varphi(x)|^2\,dx
  -\frac1{q+1}\int_{\R^{n-1}}\varphi(0,\tilde x)^{q+1}\,d\tilde x
\end{equation}
to show that $u$ is time independent.

Assume that (\ref{bound-unu}) fails.
Since $u\leq C_u$ for some $C_u>0$,
we can find $(x_k,t_k)$ such that
$|\nabla u(x_k,t_k)|\to\infty$.
Set 
$$ M(x,t):=u^{q-1}(x,t)+|\nabla u(x,t)|^{(q-1)/q}, $$
$M_k:=M(x_k,t_k)$ and $\lambda_k:=1/M_k$.
In the same way as in the proof of Theorem~\ref{thm3},
the Doubling Lemma 
guarantees that 
we may assume 
\begin{equation} \label{M-bound}
 M(x,t)\leq 2M_k\quad\hbox{whenever}\quad
   |x-x_k|+\sqrt{|t-t_k|}\leq \frac{k}{M_k}
\end{equation}
and then the Liouville theorem for the linear heat equation
\cite[Theorem 1]{Eid}
implies that we may assume $c_k:=x_{k,1}M_k\to c_\infty\in[0,\infty)$.

Assumption (\ref{sym}) guarantees
$u(x_1,\tilde x,t)=v(x_1,r,t)$, where $r=|\tilde x|$.
Passing to a subsequence we may assume
$\rho_k:=r_k/\lambda_k\to \rho_\infty\in[0,\infty]$, where $r_k=|\tilde x_k|$.
If $\rho_\infty=\infty$ then we set
$$ w_k(y,\rho,s):=\lambda_k^{1/(q-1)}v(\lambda_k y,
     r_k+\lambda_k \rho,    
     t_k+\lambda_k^2 s),
   \quad y\geq0,\ \rho\geq-\rho_k,\ s\in\R.$$
Then 
$$ \begin{aligned}
 & w_k^{q-1}(c_k,0,0)+|\nabla w_k(c_k,0,0)|^{(q-1)/q} =1, \\
 & w_k^{q-1}+|\nabla w_k| \leq 2,\quad
   \hbox{whenever}\quad \sqrt{(y-c_k)^2+\rho^2}+\sqrt{|s|}\leq k
\end{aligned}$$
and $w_k$ satisfy the equation
$$ w_s-w_{yy}-w_{\rho\rho}=\frac{n-2}{\rho_k+\rho}w_\rho,
\qquad y>0,\ \rho>-\rho_k \ \ s\in\R,$$
and the boundary condition $ w_y=-w^q$ for $y=0$.
Consequently,
a subsequence of $\{w_k\}$ converges to
a nonnegative nontrivial solution of (\ref{eq-ubc}) with $n=2$
which contradicts Theorem~\ref{thm3}.
Hence $\rho_\infty<\infty$.
Set
$$ v_k(y,s):=\lambda_k^{1/(q-1)}u(\lambda_k y,
         t_k+\lambda_k^2 s),
   \quad y\in\R^n_+,\ s\in\R,$$
fix $\tilde y\in\R^{n-1}$ with $|\tilde y|=1$
and set $y_k=(c_k,\rho_k\tilde y)$.
Then $v_k$ are solutions of (\ref{eq-ubc}) 
satisfying  (\ref{sym}), 
$v_k^{(q-1)}(y_k,0)+|\nabla v_k(y_k,0)|^{(q-1)/q}=1$
and the bound (\ref{M-bound}) guarantees that
a suitable subsequence of $\{v_k\}$ converges
to a positive bounded solution $v$
of (\ref{eq-ubc}) satisfying (\ref{sym})
and having bounded spatial derivatives. Replacing $u$ by $v$
we obtain (\ref{bound-unu}). 

Next we use doubling and scaling arguments together with Theorem~\ref{thm3}
in order to show
\begin{equation} \label{decay-tilde}
u(0,\tilde x,t)\leq C|\tilde x|^{-1/(q-1)}
\quad\hbox{for all }\ \tilde x\in\R^{n-1}.
\end{equation}
Notice that the monotonicity property (\ref{monotone}) will then  guarantee
\begin{equation} \label{decay-tilde2}
u(x,t)\leq C|\tilde x|^{-1/(q-1)}
\quad\hbox{for all }\ \tilde x\in\R^n_+.
\end{equation}
Assume on the contrary that (\ref{decay-tilde}) fails.
Then there exist $\tilde x_k,t_k$ such that  
$$u(0,\tilde x_k,t_k)|\tilde x_k|^{1/(q-1)}\to\infty.$$
Due to (\ref{bound-unu}) we have $|\tilde x_k|\to\infty$.
Denote $r=|\tilde x|$, $v(x_1,r,t)=u(x_1,\tilde x,t)$, $r_k=|\tilde x_k|$
and $M(r,t)=v(0,r,t)^{q-1}$ for $(r,t)\in(0,\infty)\times\R$.
Then $M(r_k,t_k)r_k\to\infty$ so that
we may assume $M_k:=M(r_k,t_k)>2k/r_k$.
Now the Doubling Lemma
(with 
$X=[0,\infty)\times\R$, $D=(0,\infty)\times\R$ and
$\hbox{dist}=\hbox{dist}_P$)
guarantees that we may also assume
\begin{equation} \label{doubleM}
M(r,t) \leq 2M_k\qquad \hbox{whenever }\  |r-r_k|+\sqrt{|t-t_k|}\leq\frac k{M_k}.
\end{equation}
Set $\lambda_k=1/M_k$ and
$$ w_k(y,\rho,s):=\lambda_k^{1/(q-1)}v(\lambda_k y,r_k+\lambda_k\rho,t_k+\lambda_k^2s),
 \quad y\geq0,\ \rho\geq-r_k/\lambda_k,\ s\in\R.$$
Then  $w(0,0,0)=1$  and
(\ref{doubleM}), (\ref{monotone}) guarantee  
$w_k\leq 2^{1/(q-1)}$ whenever $|\rho|+\sqrt{|s|}\leq k$.
In addition $w=w_k$ is a positive solution of the equation
$$ w_s-w_{yy}-w_{\rho\rho}=\frac{n-2}{r_k/\lambda_k+\rho}w_\rho$$
complemented by the boundary condition $ w_y=-w^q$ for $y=0$.
Since $r_k/\lambda_k\to\infty$, it is easy to pass
to the limit (in the weak formulation of the problem) 
to obtain a positive bounded solution of the problem
$w_t-\Delta w=0$ in $\R^2_+\times\R$, $w_\nu=w^q$ on $\partial\R^2_+\times\R$,
which contradicts Theorem~\ref{thm3}.
Consequently, (\ref{decay-tilde}) and (\ref{decay-tilde2}) are true.

To prove the decay of $u$ with respect to $x_1$ we use the representation
formula
\begin{equation} \label{AA}
\left. \begin{aligned}
u(x,&t)=\int_{\R^n_+}G(x,y,t-T)u(y,T)\,dy \\
  &\quad
  +\int_T^t\int_{\R^{n-1}}\partial_{y_1}G(x,(0,\tilde y),t-s)u((0,\tilde y),s)\,d\tilde y\,ds
  =: A_1+A_2,
\end{aligned}\ \right\} 
\end{equation}
for $x_1>0$ and $t>T$,  where 
$$ G(x,y,t)=\frac1{(4\pi t)^{n/2}}\bigl(e^{-|x-y|^2/4t}-e^{-|x'-y|^2/4t}\bigr),
\qquad x':=(-x_1,\tilde x).$$
Notice that 
$$\begin{aligned}
0 \leq G(x,y,t) &\leq Ct^{-n/2}e^{-|x-y|^2/4t}, \\
0 \leq \partial_{y_1}G(x,(0,\tilde y),t) &\leq Cx_1t^{-n/2-1}e^{-|x-y|^2/4t}.
\end{aligned} $$
Introducing the new variable $z=(x-y)/2\sqrt{t-T}$ in $A_1$ we have
$$ A_1\leq C\int_{\{z:z_1\leq x_1/2\sqrt{t-T}\}} e^{-|z|^2}u(x-2z\sqrt{t-T},T)\,dz \to 0
 \quad\hbox{as }\ T\to-\infty,$$
due to the Lebesgue dominated convergence theorem and the pointwise convergence
$u(x-2z\sqrt{t-T},T)\to0$ for $\tilde z\ne0$ (which follows from (\ref{decay-tilde2})).
Using estimate \eqref{decay-tilde2} we also have
$$ A_2\leq C\int_T^t x_1(t-s)^{-3/2}e^{-x_1^2/4(t-s)} I(\tilde x,t-s)\,ds,$$
where
$$ I(\tilde x,t):=\int_{\R^{n-1}}(4t)^{-(n-1)/2}e^{-|\tilde x-\tilde y|^2/4t} 
                                             |\tilde y|^{-1/(q-1)}\,d\tilde y.$$
Introducing the variable $\tilde z=(\tilde x-\tilde y)/2\sqrt{t}$ we obtain
$$ \begin{aligned}
 &I(\tilde x,t)
 = \int_{\R^{n-1}}e^{-|\tilde z|^2}|\tilde x-2\sqrt{t}\tilde z|^{-1/(q-1)}\,d\tilde z \\
 &\ \leq \int_{|\tilde x-2\sqrt{t}\tilde z|>\sqrt{t}}
     e^{-|\tilde z|^2}t^{-1/2(q-1)}\,d\tilde z
 + \int_{|\tilde x-2\sqrt{t}\tilde z|<\sqrt{t}}
   |\tilde x-2\sqrt{t}\tilde z|^{-1/(q-1)}\,d\tilde z \\
 &\ \leq C_1t^{-1/2(q-1)}+C_2\int_0^{1/2}(r\sqrt{t})^{-1/(q-1)}\cdot r^{n-2}\,dr \\
 &\ \leq C_3t^{-1/2(q-1)}.
\end{aligned}$$ 
Consequently,
$$ A_2\leq C\int_{-\infty}^t x_1(t-s)^{-3/2}(t-s)^{-1/2(q-1)}e^{-x_1^2/4(t-s)}\,ds $$
and introducing the new variable $\tau$ satisfying $t-\tau=(t-s)/x_1^2$
we obtain
$$ A_2\leq Cx_1^{-1/(q-1)}\int_{-\infty}^t (t-\tau)^{-3/2-1/2(q-1)}e^{-1/4(t-\tau)}\,d\tau
  = C(q)x_1^{-1/(q-1)}.$$
Since the last estimate of $A_2$ does not depend on $T$ and $A_1\to0$ as $T\to-\infty$,
(\ref{AA}) and (\ref{decay-tilde2}) imply
\begin{equation} \label{decay-u}
  u(x,t)\leq C|x|^{-1/(q-1)}\quad \hbox{for all }\ (x,t)\in\R^n_+\times\R.
\end{equation}
Next we use doubling and scaling arguments again to prove the estimate
\begin{equation} \label{decay-nabla}
 |\nabla u(x,t)|\leq C|x|^{-q/(q-1)}\quad \hbox{for all }\ (x,t)\in\R^n_+\times\R.
\end{equation}
Assume on the contrary that there exist $x_k,t_k$ such that
$$|\nabla u(x_k,t_k)|\cdot|x_k|^{q/(q-1)}\to\infty.$$
Due to (\ref{bound-unu}) we have $|x_k|\to\infty$.
Set 
$$ M(x,t):=u(x,t)^{q-1}+|\nabla u(x,t)|^{(q-1)/q}. $$
Then without loss of generality we may assume
$M_k:=M(x_k,t_k)>2k/|x_k|$.
The Doubling Lemma 
(with $X=\overline{\R^n_+}\times\R$,
$D=(\overline{\R^n_+}\setminus\{0\})\times\R$
and 
$\hbox{dist}=\hbox{dist}_P$)
shows that we may assume 
$$ M(x,t)\leq 2M_k\quad\hbox{whenever }\ |x-x_k|+\sqrt{|t-t_k|}\leq\frac k{M_k}.$$
Finally, we may also assume that
$c_k:=x_{k,1}M_k\to c_\infty\in[0,\infty]$, where 
$x_{k,1}$ denotes the first component of $x_k$.
Set $\lambda_k=1/M_k$ and 
$$ v_k(y,s):=\lambda_k^{1/(q-1)}u(x_k+\lambda_ky,t_k+\lambda_k^2s),\quad
  y\in\R^n,\ y_1\geq-c_k,\ s\in\R.$$ 
Then $v_k$ solves the linear heat equation and satisfies the boundary condition
$v_\nu=v^q$.
Since (\ref{decay-u}) implies $u(x,t)^{q-1}\leq C/|x|$, we have
$|\nabla u(x_k,t_k)|^{(q-1)/q}>M_k/2$ for $k$ large enough, hence
$|\nabla v_k(0,0)|>2^{-q/(q-1)}$.
On the other hand, 
$$v_k(0,0)\leq C(M_k|x_k|)^{-1/(q-1)}\to0$$
and 
$$v_k^{q-1}+|\nabla v_k|^{(q-1)/q}\leq 2\quad\hbox{for }\ |y|+\sqrt{|s|}\leq k, \
\ y_1\geq -x_{k,1}/\lambda_k.$$
If $c_\infty=\infty$ then a suitable subsequence of $v_k$
converges to the nonnegative solution $v$ of the linear heat equation in $\R^n\times\R$,
$v(0,0)=0$ and $\nabla v(0,0)\ne0$ which contradicts \cite[Theorem 1]{Eid}.
If $c_\infty<\infty$ then
a subsequence of $v_k$
converges to the nonnegative solution $v$ of the linear heat equation in 
$\{y\in\R^n:y_1>-c_\infty\}\times\R$ satisfying the boundary condition
$v_\nu=v^q$ and $v(0,0)=0$, $\nabla v(0,0)\ne0$, which yields a contradiction again.
Consequently, (\ref{decay-nabla}) is true.

Estimates (\ref{decay-u}) and (\ref{decay-nabla}) guarantee
that the energy  $E(u(\cdot,t))$ is well defined
and that we can use the same arguments as in 
the proof of Theorem~\ref{thmUa}
to show $E(u(\cdot,t))\equiv0$. Consequently, $u$ is time-independent
which contradicts the elliptic Liouville theorem in \cite{Hu}.
\qed

\vskip5mm
{\it Proof of Theorem~\ref{thm3b}.\/}
Set
$$ M(t):=\max_{x\in\overline\Omega,\,\tau\in[T/4,t]}
 \bigl( u(x,\tau)^{q-1}+|\nabla u(x,\tau)|^{(q-1)/q}\bigr),
\quad t\in[T/4,T).$$
We will prove $M(t)\sqrt{(T-t)}\leq C=C(u)$ for $t\in(T/2,T)$.

Assume on the contrary that there exist $t_k\in(T/2,T)$ 
such that 
$$M(t_k)\sqrt{T-t_k}\to\infty.$$
We may assume $M(t_k)\sqrt{T-t_k}>2k$ and $M(t_k)>2M(T/2)$.
Using the Doubling Lemma 
(with $X=[T/4,T]$, $D=[T/4,T)$ and 
$\hbox{dist}(t,\tilde t)=\sqrt{|t-\tilde t|}$)
we find $\tilde t_k\in[T/4,T)$ such that
\begin{equation} \label{M1}
 M_k:=M(\tilde t_k)\geq M(t_k), \qquad
 M_k\sqrt{T-\tilde t_k}>2k 
\end{equation}
and
\begin{equation} \label{M2}
 M(t)\leq 2M_k \quad\hbox{for all}\quad
  t\in[T/4,T),\ \sqrt{|t-\tilde t_k|}\leq\frac{k}{M_k}.
\end{equation}
In fact, the monotonicity of $M$ and (\ref{M2}) guarantee
\begin{equation} \label{M3}
 M(t)\leq 2M_k \quad\hbox{for all}\quad
  t\in\Bigl[T/4,\tilde t_k+\frac{k^2}{M_k^2}\Bigr).
\end{equation}
Inequalities $M_k\geq M(t_k)>2M(T/2)$ guarantee 
$\tilde t_k>T/2$.
Fix 
$x_k\in\overline\Omega$ and $\tau_k\in[T/4,\tilde t_k]$ such that 
$$ M_k=u(x_k,\tau_k)^{q-1}+|\nabla u(x_k,\tau_k)|^{(q-1)/q}$$
and notice that $M(\tau_k)=M_k$.
Next we distinguish two cases:
$$ \begin{aligned}
&\hbox{(i) } u(x_k,\tau_k)^{q-1}>\frac12 M_k, \\
&\hbox{(ii) } u(x_k,\tau_k)^{q-1}\leq \frac12 M_k.
\end{aligned}$$

Case (i): Since $u(x_k,\tau_k)^{q-1}>M_k/2>M(T/2)$, we have 
$u(x_k,\tau_k)>\max\{u(x,t):x\in\overline\Omega,\ t\in[T/4,T/2]\}$
and the maximum principle guarantees that there exist $\hat x_k\in\partial\Omega$
and $\hat\tau_k\in(T/2,\tau_k]\subset[T/4,\tilde t_k]$ such that 
$$ u(\hat x_k,\hat\tau_k)
     =\max\{u(x,t):x\in\overline\Omega,\ t\in[T/4,\tau_k]\}
     \geq u(x_k,\tau_k).$$
Consequently, $u^{q-1}(\hat x_k,\hat\tau_k)>M_k/2$.
Set 
$$\lambda_k:=1/M_k,\quad
\Omega_k:=\{y\in\R^n:
 \hat x_k+\lambda_k R_k y\in\Omega\},
$$ 
where $R_k$ is a rotation operator such that
$(-1,0,0,\dots,0)$ is the exterior normal vector of $\partial\Omega_k$ at $0$.
Given $y\in\overline\Omega_k$ and 
$s\in I_k:=\{s:\hat\tau_k+\lambda_k^2s\in[T/4,T)\}$,
set also
\begin{equation} \label{vk}
 v_k(y,s):=\lambda_k^{1/(q-1)}u(\hat x_k+\lambda_k R_k y,\hat\tau_k+\lambda_k^2s).
\end{equation}
Then $v_k$ solve the equation and the boundary condition
in (\ref{eq-nbc}) in $\Omega_k\times I_k$ and 
on $\partial\Omega_k\times I_k$, respectively,   
$v_k^{q-1}(0,0)>1/2$, $v^{q-1}(y,s)+|\nabla v(y,s)|^{(q-1)/q}\leq 2$
for all $(y,s)\in \overline\Omega_k\times I_k$ satisfying $|s|\leq k^2$.
The arguments in \cite{Hu96} guarantee that a subsequence of $\{v_k\}$
converges in $C^1_{loc}$ to a positive entire solution of (\ref{eq-ubc})
which contradicts Theorem~\ref{thm3}.

Case (ii):
Denote $d_k:=\hbox{dist}(x_k,\partial\Omega)$
and choose $\tilde x_k\in\partial\Omega$ such that $|\tilde x_k-x_k|=d_k$.
Set also $\hat\tau_k:=\tau_k$
and $c_k:=d_kM_k$. We may assume $c_k\to c_\infty\in[0,\infty]$.

If $c_\infty<\infty$  
then we set $\hat x_k:=\tilde x_k$, 
$y_k:=(c_k,0,0,\dots,0)$ and
define $\lambda_k$,$\Omega_k$,\allowbreak$R_k$,$I_k$,$v_k$ as in Case (i).
Notice that  $R_ky_k=\tilde x_k-\hat x_k$.
Similarly as in Case (i),
$v_k$ solve the equation and the boundary condition
in (\ref{eq-nbc}) in $\Omega_k\times I_k$ and 
on $\partial\Omega_k\times I_k$, respectively,   
$v_k^{q-1}(y_k,0)+|\nabla v_k(y_k,0)|^{(q-1)/q}=1$, 
$v^{q-1}(y,s)+|\nabla v(y,s)|^{(q-1)/q}\leq 2$
for all $(y,s)\in \overline\Omega_k\times I_k$ satisfying $|s|\leq k^2$
and the arguments in \cite{Hu96} guarantee that a subsequence of $\{v_k\}$
converges in $C^1_{loc}$ to a positive entire solution of (\ref{eq-ubc})
which contradicts Theorem~\ref{thm3}.

If $c_\infty=\infty$
then we set $\hat x_k:=x_k$, define $R_k$ as the identity
and $\lambda_k,\Omega_k,I_k,v_k$ as in Case (i).
Now a subsequence of $\{v_k\}$ converges 
in $C^1_{loc}$ to a nonnegative bounded solution 
of the linear heat equation in $\R^n\times\R$.
Since  
$$|\nabla u(\hat x_k,\hat\tau_k)|^{(q-1)/q}
=|\nabla u(x_k,\tau_k)|^{(q-1)/q}\geq M_k/2$$
we have $|\nabla v_k(0,0)|\geq 1/2$, hence $v$ is nonconstant
which contradicts the Liouville theorem for the linear heat equation
\cite[Theorem 1]{Eid}.
\qed


\end{document}